\documentclass[12pt]{article}
\usepackage{amsmath,amsthm,amsfonts,amssymb,amscd}
\usepackage[latin1]{inputenc}

\renewcommand{\thefootnote}

\def\_#1{{\lower 0.7ex\hbox{}}_{#1}}

\begin{document}

\theoremstyle{plain}
\newtheorem{thm}{\sc Theorem}[section]
\newtheorem{lem}[thm]{\sc Lemma}

\newtheorem*{coro}{\sc Corollary}
\newtheorem{prop}[thm]{\sc Proposition}

\theoremstyle{remark}
\newtheorem{Case}{\bf Case}[section]
\newtheorem{rem}[thm]{\sc Remark}
\newtheorem*{exem}{\sc Example}

\theoremstyle{definition}
\newtheorem{Def}[thm]{\bf Defini\c c\~ao}
\newtheorem*{pf}{\sc Proof}
\newtheorem{obs}[thm]{\bf Observa\c c\~ao}

\numberwithin{equation}{section}

\def\argmin{\operatorname{argmin}}
\def\grad{\operatorname{grad}}
\def\tni{\operatorname{int}}
\def\ir{\operatorname{ir}}
\newcommand{\ve}{{\varepsilon}}
\newcommand{\be}{{\beta}}
\newcommand{\la}{{\lambda}}
\newcommand{\teta}{{\theta}}
\newcommand{\Om}{{\Omega}}
\newcommand{\La}{{\Lambda}}
\newcommand{\ga}{{\gamma}}
\newcommand{\Ga}{{\Gamma}}
\newcommand{\po}{{\partial}}
\newcommand{\ro}{{\rho}}
\newcommand{\ov}{\overline}
\newcommand{\re}{{\mathbb{R}}}
\newcommand{\OI}{{\mathcal O}}
\newcommand{\te}{{\theta}}
\newcommand{\lgg}{{\left\langle\right.}}
\newcommand{\rg}{{\left.\right\rangle}}
\newcommand{\lV}{{\left\Vert \right.}}
\newcommand{\rV}{{\left.\right\Vert}}
\newcommand{\na}{{\mathbb N}}
\newcommand{\al}{{\alpha}}
\newcommand{\lv}{{\left\vert \right.}}
\newcommand{\rv}{{\left.\right\vert}}
\newcommand{\vr}{{\varphi}}
\newcommand{\ck}{{\mathcal C}}
\newcommand{\om}{{\omega}}

\begin{title}
{\bf Complete hypersurfaces in Euclidean spaces with strong finite
total curvature}
\end{title}

\vskip .5in

\author{
Manfredo do Carmo\,\, and\,\, Maria Fernanda Elbert}

\date{}
\maketitle

 \footnote{{\bf Key words and sentences:} Complete
hypersurface, total curvature, topological properties,
 Gauss-Kronecker curvature \par
\par {\bf Both authors are partially supported by  CNPq and Faperj}.}

\begin{abstract}
We prove that strong finite total curvature (see definition in Section 2)
complete hypersurfaces of $(n+1)$-euclidean space are proper and
diffeomorphic to a compact manifold minus finitely many points.
With an additional condition, we also prove that the Gauss map of
such hypersurfaces extends continuously to the punctures. This is
related to results of White \cite{W} and
and
Müller-$\check{\text{S}}$verák  \cite{MS}. Further properties of
these hypersurfaces are presented, including a gap theorem for the
total curvature.
\end{abstract}

\vskip .5in

\section{Introduction}

Let $\phi\colon M^n\to\re^{n+1}$ be a hypersurface of the euclidean
space $\re^{n+1}$. We assume that $M^n=M$ is orientable and we fix an
orientation for $M$. Let $g\colon M\to S^n_1\subset \re^{n+1}$ be
the Gauss map in the given orientation, where $S_1^n$ is the unit
$n$-sphere. Recall that the linear operator $A\colon T_p M\to T_p
M$, $p\in M$, associated to the second fundamental form, is given
by
$$
\langle A(X),Y\rangle=-\langle\ov\nabla_X N,Y\rangle,\quad X,Y\in
T_pM,
$$
where $\ov\nabla$ is the covariant derivative of the ambient space
and $N$ is the unit normal vector in the given orientation. The
map $A=-dg$ is self-adjoint and its eigenvalues are the principal
curvatures $k_1,k_2,\dots,k_n$.

 We say that the total
curvature of the immersion is finite if $\int_M|A|^n\, dM<\infty$, where $|A|=\big(\sum_ik_i^2\big)^{1/2}$, i.e., if $|A|$ belongs to the space 
$L^n(M)$. If $\phi\colon M^n\to\re^{n+1}$ is a complete minimal hypersurface
with finite total curvature then $M$ is (equivalent to) a compact
manifold $\ov M$ minus finitely many points and the Gauss map
extends to the punctures. This was proved by Osserman \cite{O} for
$n=2$ (the equivalence here is conformal and the Gauss map extends
to a (anti) holomorphic map $\ov g\colon\ov M^2\to S^2_1$; the
conformal equivalence had already been proved by Huber \cite{Hu}).
For an arbitrary $n$, this was proved by Anderson \cite{A} (here
the equivalence is a diffeomorphism and the Gauss map extends
smoothly).

When $\phi$ is not necessarily minimal and $n=2$, the above result,
with the additional hypothesis that the Gauss curvature does not
change sign at the ends, was shown to be surprisingly true by B.
White \cite{W}. The subject was taken up again by Müller-$\check{\text{S}}$verák  \cite{MS} who answered a question of
\cite{W} and obtained further information on the conformal
behaviour of the ends.

The results of White \cite{W} and  Müller-$\check{\text{S}}$verák \cite{MS} start from the fact that,
since $\int_{M^2}|A|^2\, dM\ge2 \int_{M^2}|K|\, dM$, finite total
curvature for $n=2$ implies, by Huber's theorem, that $M$ is homeormorphic
to a compact surface minus finitely many points.
 For an arbitrary dimension, any generalization of Huber's theorem should require stronger assumptions (see \cite{CH1} and \cite{CH2} for a discussion on the theme). Thus, for a generalization of \cite{W} and \cite{MS}  for $n \ge 3$ a further condition might be necessary to account for the lack of an appropriate generalized Huber
 theorem.

  Here, we assume the hypothesis of {\it strong finite total curvature}, that is, we assume that $|A|$ belongs to $W^{1,q}_{s}$, a special Weighted Sobolev space (see Section 2 for precise definitions). We point out that the spaces $W^{1,q}_{s}$ were used in a seminal work of R. Bartnik \cite{B} for establishing a decay condition on the metric of an n-manifold, $n\geq 3$,  in order to prove that the ADM-mass is well-defined. Following the ideas of \cite{B}, a lot of related papers also uses the norm of $W^{1,q}_{s}$ to express decay assumptions (see for instance \cite{LF}, \cite{DM}, \cite{OW}).

  \

We have proved the following results.

\begin{thm}\label{thm1.1}
Let $x\colon M^n\to\re^{n+1}$, $n\ge3$, be an orientable, complete
hypersurface with strong finite total curvature. Then:

\begin{itemize}\itemsep=-1pt
\item[{\rm i)}] The immersion $\phi$ is proper. \item[{\rm ii)}] $M$
is diffeomorphic to a compact manifold $\ov M$ minus a finite
number of points $q_1, \dots q_k$.
\end{itemize}

Assume, in addition, that the Gauss-Kronecker curvature
$H_n=k_1k_2\dots k_n$ of $M$ does not change sign in punctured
neighbourhoods of the $q_i$'s.  Then:

\begin{itemize}\itemsep=-1pt
\item[{\rm iii)}] The  Gauss map $ g\colon M^n \to S^n_1$ extends
continuously to the points $q_i$.
\end{itemize}
\end{thm}

\begin{thm}\label{thm1.2}
Let  $x\colon M^n\to\re^{n+1}$, $n \ge 3$, be an orientable
complete hypersurface with strong finite total curvature. Assume
that the set $N$ of critical values of the Gauss map $g$ is a
finite union of submanifolds of $S_1^n$ with codimension $\ge 3$.
Then:
\begin{itemize}\itemsep=-1pt
\item[\rm i)] The extended Gauss map $\bar g\colon \ov M \to
S_1^n$ is a homeomorphism. \item[\rm ii)] If, in addition, $n$ is
even, $M$ has exactly two ends.
\end{itemize}
\end{thm}

\noindent{\bf Remark}.The condition on $N$ can be replaced by a
weaker condition on the Hausdorff dimension of $N$ and the rank of
$g$ (See \cite{MO}, Theorems B and C and Remark 6.7).

\

 It follows from Theorem \ref{thm1.1} that
there is a computable lower bound for the total curvature of the
non-planar hypersurfaces of the set $C^n$ defined in the statement
below.

\begin{thm}\label{thm1.3} {\rm (The Gap Theorem)} Let $C^n$ be the set of strong finite total curvature complete orientable hypersurfaces $\phi\colon M^n \to \re^{n+1}$, $n \ge 3$, such that $H_n$ does not change sign in $M$. Then either $x(M^n)$ is a hyperplane, or
$$
\int_M |A|^n\,dM > 2 \sqrt{n!} \,
(\sqrt\pi)^{n+1}\big/\Ga((n+1)/2)),
$$
where $\Ga$ is the gamma function.
\end{thm}

\medskip

\noindent{\bf Remark}. For the Gap Theorem it is not enough to requiring that $H_n$ does not change sign at the ends of the
hypersurface. This should hold on the whole $M$. Consider the rotation hypersurfaces in $\re^{n+1}$
generated by the smooth curve $x^{n+1} = \ve\,e^{-1/x_1^2}$, $\ve
> 0$, $(x_1,\dots,x_n, x_{n+1}) \in \re^{n+1}$. It is easily
checked that, for all $\ve$, this hypersurface has strong finite
total curvature  and $H_n$ does not change sign at the (unique)
end of the hypersurface. However, as $\ve$ approaches zero, these
hypersurfaces approach a hyperplane, and the lower bound for the
total curvature of the family is zero.

\vskip .1in

The paper is organized as follows. In Section 2, we discuss
(Proposition 2.2) the rate of decay at infinity of the second
fundamental form of a hypersurface under the hypothesis of  strong
finite total curvature. In Section 3, we show that each end of
such a hypersurface has a unique ``tangent plane at infinity''
(see the definition before Proposition 3.4) and in Section 4, we
prove Theorems 1.1, 1.2 and the Gap Theorem.

\section{The rate of decay of the second fundamental form}

In the rest of this paper, we will be using the following notation
for an immersion \linebreak $\phi\colon M^n\to\re^{n+1}$:
\begin{align*}
\ro &= \text{intrinsic distance in $M$}\\
d &= \text{distance in $\re^{n+1}$;\,\, $0=$ origin of $\re^{n+1}$}\\
D_p(R) &=\{x\in M;\,\, \ro(x,p)<R\}\\
D_p(R,S) &=\{x\in M;\,\, R<\ro(x,p)<S\}\\
B(R) &=\{x\in \re^{n+1}; \,\, d(x,0)<R\}; \,\, S(R)=\partial B(R)\\
A(R,S) &=\{x\in\re^{n+1};\,\, R<d(x,0)< S\}.
\end{align*}

\vskip .1in

 Without loss of generality, we assume that $0\in \phi(M)$ and we choose a point  $0\in M$ such that $0=\phi(0)$. For $x\in M$, $\ro_0 (x)$ will denote the intrinsic distance in $M$ from $x$ to $0$. Now, we set the notation for the norms (see \cite[(1.2)]{B}) that will be used in the definition of strong finite total curvature.

Let $\Omega\subset M$.  Given any $q>0$, we define the {\it weighted space } $L^q _s (\Omega)$ of all measurable functions of finite norm
$$
||u||_{L^q_s (\Omega)}=\left ( \int_\Omega |u|^q |\ro_0|^{-qs-n} \;dM\right)^{1/q}\hspace{-15pt}.
$$

For positive integers $k$, we introduce the {\it weighted Sobolev space} $W^{k,q} _s(\Omega)$ of all measurable functions of finite norm
$$
||u||_{W^{k,q} _s}=\sum ^k _{i=0}||\nabla u||_{L^q _{s-i} (\Omega)}.
$$

We say that the immersion {\bf \it has strong finite total curvature} if
$$|A|\in W^{1,q} _{-1}(M)\;\;\mbox{for}\; q>n,$$

 that is, if
$$
\left ( \int_\Omega |A|^q |\ro_0|^{q-n} \;dM\right)^{1/q}+\left ( \int_\Omega |\nabla |A||^q |\ro_0|^{2q-n} \;dM\right)^{1/q}<\infty\;\;\mbox{for}\; q>n.
$$
\vskip .1in

We remark that the function $\ro_0$ used above to define these norms could be replaced by the distance with respect to any point $p\in M$. Here, we fixed the point $0=\phi(0)$. Then, in this paper, when we say that the immersion has strong finite total curvature we are implicitly assuming w.l.g. that $0\in\phi(M)$. We also remark that the weights used to define the norm $||.||_{W^{k,q} _s}$ makes it invariant by dilations (see the proof of Proposition \ref{prop2.1}).

\

 The following Lemma will be repeatedly used in this and in the next
section.

\begin{lem}\label{lem2.1}
Let $D\subset \re^{n+1}$ be a bounded domain with smooth boundary
$\po D$. Let $(W_i)$ be a sequence of connected n-manifolds and
let $\phi\colon W_i\to \re^{n+1}$ be immersions such that $x(\po
W_i)\cap D=\emptyset$ and $x(W_i) \cap D = M_i$ is connected and
nonvoid. Assume that there exists a constant $C>0$ such that
${\displaystyle\sup_{x\in M_i}|A_i(x)|^2<C}$ and that there exists a sequence of
points $(x_i)$, $x_i\in M_i$, with a limit point $x_0\in D$. Then:
\begin{itemize}\itemsep=-1pt
\item[\rm i)] A subsequence of $(M_i)$ converges $C^1$ uniformly
on the compact parts (see the definition below) to   a union of
hypersurfaces $M_\infty \subset D$.
\item[\rm ii)] If, in
addition, $\left ( \int_\Omega |A|^q \alpha_i \;dM\right)^{1/q}+\left ( \int_\Omega |\nabla |A||^q \beta_i \;dM\right)^{1/q}\to 0$,  for sequences $(\alpha_i)_i$ and $(\beta_i)_i$ of continuous functions such that ${\displaystyle\inf_{x\in
M_i}\{\alpha_i,\beta_i\}\geq \kappa}>0$. Then a subsequence of $|A_i|$ converges to zero everywhere and
$M_\infty$ is a union of hyperplanes.
\end{itemize}
\end{lem}

By $C^1$ convergence to $M_\infty$ we mean that for any $m\in
M_\infty$ and each tangent plane $T_m M_\infty$ there exists an
euclidean ball $B_m$ around $m$ so that, for $i$ large, the image
by $\phi$ of some connected component of $\phi^{-1}(B_m\bigcap M_i)$ can
be graphed over $T_m M_\infty$ by a function $g_i^m$ and the
sequence $g_i^m$ converges $C^1$ to the graph $g_\infty$ of
$M_\infty$ over the chosen plane $T_m M_\infty$.
\begin{pf}
From the uniform bound of the curvature $|A_i|^2$, we conclude the
existence of a number $\delta>0$ such that for each $p_i\in M_i$
and for each tangent space $T_{p_i}M_i$, $M_i$ can be graphed by a
function $f_i^{p_i}$ over a disk $U_\delta(p_i)\subset
T_{p_i}M_i$, of radius $\delta$ and center $p_i$ in $T_{p_i}M_i$,
and that such functions have a uniform $C^1$ bound (independent of
$p_i$ and $i$). We want to show that we also have a uniform $C^2$
bound.
\end{pf}

Let $q$ be a point in the part of $M_i$ that is a graph over
$U_{\delta}(p_i)$ and let $v \in T_qM_i$\,. Consider the plane
$P_q$ that contains the normal vector $N_i(q)$ and $v$ and take
the curve $C_i = P_q \cap M_i$\,. Parametrize $C_i$ by $c_i(t)$
with $c_i(0)=q$, project it down to $T_{p_i}M_i$ parallely to the
normal at $p_i$\,. Let $\tilde c_i(t)$ be this projection; then,
$c_i(t) = \big(\tilde c_i(t), f_i^{p_i}(\tilde c_i(t))\big)$ and
the normal curvature of $M_i$ in $q$ along $v$ is
\begin{equation}
k_v^i(q) =
\big(f_i^{p_i}\big)''(0)\big/\big(1+\big[\big(f_i^{p_i}\big)'(0)\big]^2\big)^{3/2},\label{extra}
\end{equation}
where, e.g., $(f_i^p)'(t)$ means the derivative in $t$ of
$f_i^{p_i}(\tilde c_i(t)) = f_i^{p_i}(t)$. It follows that we have
a uniform estimate for second derivatives in any direction $v$. By
a standard procedure (see e.g. \cite{CS} p. 280), this implies a
uniform $C^2$-bound on $f_i^{p_i}$.  Now, consider the sequence
$(x_i)$ with a limit point $x_0$, and let $\tau_i$ be the
translation that takes $x_i$ to $x_0$. The unit normals of
$\tau_i(M_i)$ at $x_0$ have a convergent subsequence, hence a
subsequence of the tangent planes $T_{x_0}(\tau_iM_i)$ converges
to a plane $P$ containing $x_0$. For $i$ large, the parts of $M_i$
that were graphs over $U_\delta(x_i)$ are now graphs over
$U_{\delta/2}(x_0)\subset  P$; we will denote the corresponding
functions by $g_i^{x_0}$. By the bounds on the derivatives that we
have obtained, the functions $g_i^{x_0}$ and their first and
second derivatives are uniformly bounded, say, $|g_i^{x_0}|\_{2;
U_{\delta/2}(x_0)} <C_1$. By the theorem of Arzel\'a-Ascoli and
the fact that, in this case, limits and derivatives are
interchangeable,  a subsequence of $g_i^{x_0}$ converges $C^1$,
uniformly in the compact subdomains of $U_{\delta/2}(x_0)$ to a
function $g_\infty^{x_0}$.

Notice that we have obtained a subsequence of $(M_i)$ with the
property that those parts of $M_i$ that are graphs around the
points $x_i$, converge to a hypersurface passing through $x_0$. We
will express this fact by saying that $(M_i)$ {\it has a
subsequence that converges locally at\/} $x_0$.

To complete the proof of (i) of  Lemma 2.1, we need a covering
argument that runs as follows.

Let $L$ be the set of all limit points of sequences of the form
$(p_i)$, where $p_i\in M_i$, and let $M_\infty$ be the connected
component of $L$ that contains $x_0$. Let $q_1,q_2,\dots$ be a
sequence of points in $M_\infty$ that is dense in $M_\infty$. Let
$(q_1^i)$, $q^i_1\in M_i$, be a sequence that converges to $q_1$.
As we did before, we can obtain a subsequence $(M_i^1)$ of $(M_i)$
that converges locally at $q_1$. From this sequence, we can
extract a subsequence $(M_i^2)$ that converges locally at $q_1$
and $q_2$. By induction, we can find sequences $(M_i^n)$ that
converge locally at $\bigcup_i q_i$, $i=1,\dots,n$. By using the
Cantor diagonal process, we obtain a sequence $M_1^1,M_2^2,\dots$
that converges $C^1$ to $M_\infty$ and shows that $M_\infty$ is a
collection of $C^1$ hypersurfaces. Clearly $M_\infty$ has no
boundary point in the interior of $D$. Thus $M_\infty$ extends to
the boundary of $D$. Since the local convergence is uniform in
compact subsets, it follows that the convergence to $M_\infty$ is
uniform in the compact subsets of $M_\infty$. This completes the
proof of (i) of Lemma 2.1.

\medskip

Now we prove (ii) of Lemma 2.1. By (i), a subsequence of $M_i$
converges $C^1$ to a  collection of hypersurfaces, $M_\infty$.
 As in the proof of (i), given $p \in M_\infty$, we can look upon the part of $M_i$ near $p$, for large $i$,
 as a graph of a function $g_i^p$ over $U_{\delta/2}(p) \subset T_pM_\infty$. The functions $g_i^p$ converge $C^1$ to the function $g^p$ that defines $M_\infty$ near $p$.

Let $G^p_i$  be the metric of $M_i$ restricted to $g^p_i(U_{\delta/2}(p))$, $G^p_\infty$  be the metric of $M_\infty$ restricted to $g^p(U_{\delta/2}(p))$ and let $E$ be the euclidean metric in $T_pM_\infty$. Notice that since the convergence $M_i\rightarrow M_\infty$ is $C^1$,  $G^p_i$ converges to $G^p_\infty$. There exists a constant $\lambda_i>0$ such that
$$
\frac{1}{\lambda_i}E(X,X)\leq G_i^p(X,X)\leq \lambda_i E(X,X),\;\;\mbox{for all}\;\;X\in T_pM_\infty\simeq\re^n.
$$

\noindent Then $dM_i=\sqrt{det(G)} dV\geq (\frac{1}{\lambda_i})^{n/2}dV,$ where $dV$ is element of
volume of $(T_pM_\infty,E)\simeq\re^n$. We obtain

$$\begin{array}{lll}
&&\left(\int_{g_i^p(U_{\delta/2}(p))} |A|^q \alpha_i \;dM\right)^{1/q}+\left ( \int_{g_i^p(U_{\delta/2}(p))} |\nabla |A||^q \beta_i \;dM\right)^{1/q}\geq\\
\\[8pt]
  &&\hspace{1.5cm}\kappa(\frac{1}{\lambda_i})^{n/2}\left(\int_{U_{\delta/2}(p)} |A|^q \;dV\right)^{1/q}+ \kappa(\frac{1}{\lambda_i})^{(n+q)/2}\left ( \int_{U_{\delta/2}(p)} |\nabla_{\hspace{-0.09cm}\mbox{\tiny E}}\; |A|\;|_{\hspace{-0.05cm}\mbox{\tiny E}}^q  \;dV\right)^{1/q}\hspace{-0.5cm}.
  \end{array}
$$

 Since
 $$
  \left(\int_{g_i^p(U_{\delta/2}(p))} |A|^q \alpha_i \;dM\right)^{1/q}+\left ( \int_{g_i^p(U_{\delta/2}(p))} |\nabla |A||^q \beta_i \;dM\right)^{1/q}\to 0$$

 \noindent we conclude that $|A_i|\to 0$ in the usual Sobolev space $W^{1,q}(U_{\delta/2}(p)).$ Now, since $q>n$, it follows from the fact that the injection  $$W^ {1,q}(U_{\delta/2}(p))\hookrightarrow C^0(U_{\delta/2}(p),\re)$$ is compact (see, for instance, \cite{AF}, page 168) that a subsequence of $(|A_i|)_i$ (again denoted by $(|A_i|)_i$) converges to zero in $||.||_{C^0}.$

Finally, we prove that $M_\infty$ is a collection of hyperplanes
by using the fact that $|A_i|\to 0$ everywhere. Since we have not
proved that the convergence is $C^2$, this is not immediate. An
argument is as follows. Let $p\in M_\infty$ and again look at the
part of $M_i$ near $p$ as a graph of a function $g_i^p$ over
$U_{\delta/2}(p)\subset T_pM_\infty$ so that, as before, $g_i^p$
converges $C^1$ to $g^p$ that defines $M_\infty$ near $p$. Let
$q\in U_{\delta/2}(p)$ and $w\in \re^n$, $|w|=1$. Set
$r(t)=q+tw\subset U_{\delta/2}(p)$, $c_i(t)=(r(t),g_i^p(r(t)))$
and $c(t)=(r(t), g^p(r(t)))$. The fact that $|A_i|\to0$ is easily
seen to imply that $(g_i^p)''(t)\to0$ in $U_{\delta/2}(p)$ (See
\eqref{extra}).

We will prove that $M_\infty$ is a hyperplane over
$U_{\delta/2}(p)$; since $p$ is arbitrary, this will yield the
result. Since we have a bound for the second derivatives of
$g_i^p$ in $U_\delta(p)$, we can use the Dominated Convergence
Theorem and the fact that $(g_i^p)'(t)\to(g^p)'(t)$ to obtain
\begin{align*}
(g^p)'(t)-(g^p)'(0) &= \lim_{j\to\infty} \{(g_i^p)'(t)-(g_i^p)'(0)\} \\
&= \lim_{j\to \infty} \int_0^t (g_i^p)''(s)\, ds = \int_0^t
\lim_{j\to\infty} (g_i^p)''(s)\, ds = 0,
\end{align*}
Thus, $c(t)$ is a straight line and, since $w$ is arbitrary,
$M_\infty$ is a hyperplane over $U_\delta(p)$, as we asserted. This concludes the proof of Lemma 2.1. \qed

\medskip

\noindent {\bf Remark.} For future use, we observe that in the
proof that $M_\infty$ is a hyperplane we only use that the
 convergence is $C^1$, that we have a bound for the second derivatives of $g_i^p$  and that $|A_i|\to 0$ everywhere.

\medskip

The proof of the following Proposition is inspired by that of
\cite{A}, Proposition 2.2; for completeness, we present it here.
Actually, the crucial point of the proof (Lemma 2.3 below), is
also similar to the proof of Proposition 2 in Choi-Schoen \cite{CS}.

\begin{prop}\label{prop2.1}
Let $\phi\colon M^n\to\re^{n+1}$ be a complete immersion with strong finite total curvature. Then, given $\ve>0$ there
exists $R_0>0$ such that, for $r>R_0$,
$$
r^2\sup_{x\in M-D_0(r)}|A|^2(x)<\ve.
$$
\end{prop}

\vskip .1in
For the two lemmas below we use the
following notation. We denote by $h\colon X^n\to\re^{n+1}$ an
immersion into $\re^{n+1}$ of an $n$-manifold $X^n=X$ with
boundary $\po X$ such that  there exists a point $x\in X$ with
$D_x(1)\cap\po X=\emptyset$.

\begin{lem}\label{lem2.3}
There exists $\delta>0$ such that if

$$
\left(\int_{D_{x_i}(1)} |A|^q \mu \;dX\right)^{1/q}+\left ( \int_{D_{x_i}(1)} |\nabla |A||^q \nu \;dX\right)^{1/q}<\delta,
$$
for any $h\colon
X^n\to\re^{n+1}$ as above and for any pair of continuous functions $\mu,\nu\colon D_x (1)\to\re$ that satisfy ${\displaystyle\inf_{D_x (1)}\{\mu,\nu\}>c>0}$,
then
$$
\sup_{t\in[0,1]}\left[t^2\sup_{D_x(1-t)}|A_h|^2\right]\le4.
$$
Here $A_h$ is the linear map associated to the second fundamental
of $h$.
\end{lem}

\begin{pf}
Suppose the lemma is false. Then there exist a sequence $h_i\colon
X_i\to \re^{n+1}$, a sequence of points $x_i\in X_i$ with $D_{x_i}(1)\cap \po
X_i=\emptyset$ and  sequences $(\mu_i)_i, (\nu_i)_i$, with  ${\displaystyle\inf_{D_x (1)}\{\mu_i,\nu_i\}>c}$ such that
$$
\left(\left(\int_{D_{x_i}(1)} |A_i|^q \mu_i \;dX_i\right)^{1/q}+\left ( \int_{D_{x_i}(1)} |\nabla |A_i||^q \nu_i \;dX_i\right)^{1/q}\right)\to 0
$$
but
$$
\displaystyle{\sup_{t\in[0,1]}\left[t^2\sup_{D_{x_i}(1-t)}|A_i|^2\right]}>4,
$$
for all $i$, where $A_i=A_{h_i}$.

Choose $t_i\in[0,1]$ so that
$$
t_i^2\sup_{D_{x_i}(1-t_i)}|A_i|^2=\sup_{t\in[0,1]}
\left[t^2\sup_{D_{x_i}(1-t)}|A_i|^2\right]
$$
and choose $y_i\in \overline{{D_{x_i}(1-t_i)}}$ so that
$$
|A_i|^2(y_i)=\sup_{D_{x_i}(1-t_i)}|A_i|^2.
$$
By using that $D_{y_i}(t_i/2)\subset D_{x_i}(1-(t_i/2))$  we obtain
$$
\sup_{D_{y_i}(t_i/2)}|A_i|^2\le\sup_{D_{x_i}(1-(t_i/2))}|A_i|^2\le
\frac{t_i^2}{t_i^2/4}\sup_{D_{x_i}(1-t_i)}|A_i|^2,
$$
hence, by the choice of $y_i$, we have
\begin{equation}\label{eq2.1}
\sup_{D_{y_i}(t_i/2)}|A_i|^2\le 4|A_i|^2(y_i).
\end{equation}

We now rescale the metric defining $d\tilde
s_i^2=|A_i|^2(y_i)ds_i^2$, that is, $d\tilde s_i^2$ is the metric
on $X_i$ induced by $\tilde h_i=d_i\circ h_i$, where $d_i$ is the
dilation of $\re^{n+1}$ about $h_i(y_i)$ (by translation, we may
assume that $h_i(y_i)=0$) by the factor
 $|A_i|(y_i)$ . The symbol $\sim$ will indicate quantities measured with
 respect to the new metric $d\tilde s_i^2$.

By assumption, $|A_i|^2(y_i)>4/t_i^2$. Thus
$$
 \widetilde D_{y_i}(1)=D_{y_i}([|A_i|(y_i)]^{-1})\subset D_{y_i}(t_i/2)\subset D_{x_i}(1-t_i/2)\subset D_{x_i}(1).
$$
It follows that $\widetilde D_{y_i}(1)\cap\po X_i=\emptyset$. Now,
we use \eqref{eq2.1} and the fact that
$$
|\widetilde A_i|(p)=[|A_i|(y_i)]^{-1}|A_i|(p)
$$
to obtain
$$
\sup_{\widetilde D_{y_i}(1)}|\widetilde A_i|^2\le 4.
$$

Therefore, the sequence $\tilde h_i=\widetilde
D_{y_i}(1)\to\re^{n+1}$, $\tilde h_i(y_i)=0$, is a sequence of
immersions with uniformly bounded second fundamental form.

By using that $
 \widetilde D_{y_i}(1)=D_{y_i}([|A_i|(y_i)]^{-1})\subset D_{x_i}(1)$ we have

$$
\begin{array}{rcl}
&&\left(\int_{D_{x_i}(1)} |A_i|^q \mu_i \;dX_i\right)^{1/q}+\left ( \int_{D_{x_i}(1)} |\nabla |A_i||^q \nu_i \;dX_i\right)^{1/q}\geq \\[14pt]
&&{\hspace{1.8cm}}\left(\int_{D_{y_i}([|A_i|(y_i)]^{-1})} |A_i|^q \mu_i \;dX_i\right)^{1/q}+\left ( \int_{D_{y_i}([|A_i|(y_i)]^{-1})} |\nabla |A_i||^q \nu_i \;dX_i\right)^{1/q}\hspace{-0.5cm}.
\end{array}$$

Thus,  we obtain

 $$
\begin{array}{rcl}
&&\left(\left(\int_{\widetilde D_{y_i}(1)} |\widetilde A_i|^q \mu_i\; |A_i(y_i)|^{q-n} \;d{\widetilde X_i}\right)^{1/q}+\left ( \int_{\widetilde D_{y_i}(1)} |\widetilde\nabla |\widetilde A_i||^q \nu_i\;|A_i(y_i)|^{2q-n} \;d{\widetilde X_i}\right)^{1/q}\right)\leq \\[16pt]
&&{\hspace{1.8cm}}\left(\left(\int_{D_{x_i}(1)} |A_i|^q \mu_i \;d X_i\right)^{1/q}+\left ( \int_{D_{x_i}(1)} |\nabla |A_i||^q \nu_i \;dX_i\right)^{1/q}\right)\to 0.
\end{array}
$$

Since $|A_i(y_i)|>\frac{2}{t_i}\geq 2$  we can use Lemma \ref{lem2.1}, with $\alpha_i=\mu_i|A_i(y_i)|^{q-n}$, $\beta_i=\nu_i|A_i(y_i)|^{2q-n}$ and $\kappa=2 c$, to conclude that a subsequence of
$|\widetilde A_i|$ converges to zero.
 But $|\widetilde A_i|(y_i)=1$, for all $i$, hence $|\widetilde A_\infty|(y_\infty)=1$. This is a contradiction, and completes the proof of Lemma \ref{lem2.3}. \qed
\end{pf}

\begin{lem}\label{lem2.4}
Given $\ve_1>0$, there exists $\delta>0$,
such that if

$$
\left(\int_{D_{x_i}(1)} |A|^q \mu \;dX\right)^{1/q}+\left ( \int_{D_{x_i}(1)} |\nabla |A||^q \nu \;dX\right)^{1/q}<\delta,
$$
for any $h\colon
X^n\to\re^{n+1}$ as above and for any pair of continuous functions $\mu,\nu\colon D_x (1)\to\re$ that satisfy ${\displaystyle\inf_{D_x (1)}\{\mu,\nu\}>c>0}$,
then

$$
\sup_{D_x(1/2)}|A_h|^2<\ve_1.
$$
\end{lem}

\vskip .1in

\begin{pf}
Suppose the lemma is false. Then there exist a sequence $h_i\colon
X_i\to \re^{n+1}$, a sequence of points $x_i\in X_i$ with $D_{x_i}(1)\cap \po
X_i=\emptyset$ and  sequences $(\mu_i)_i, (\nu_i)_i$, with  ${\displaystyle\inf_{D_x (1)}\{\mu_i,\nu_i\}>c}$ such that
\begin{equation}
\left(\left(\int_{D_{x_i}(1)} |A_i|^q \mu_i \;dX_i\right)^{1/q}+\left ( \int_{D_{x_i}(1)} |\nabla |A_i||^q \nu_i \;dX_i\right)^{1/q}\right)\to 0
\label{eq2.2}\end{equation}
but
\begin{equation}\label{eq2.3}
\sup_{D_{x_i}(1/2)}|A_i|^2\ge K^2,
\end{equation}
for some constant $K$.

By Lemma \ref{lem2.3} (with $t=1/2$), we have, for $i$
sufficiently large,
$$
\sup_{D_{x_i}(1/2)}|A_i|^2\le 16.
$$
By  \eqref{eq2.2} and Lemma \ref{lem2.1} , a subsequence of
$|A_i|$ converges to zero. This is a contradiction to
\eqref{eq2.3} and proves Lemma~\ref{lem2.4}. \qed
\end{pf}

\bigskip

\noindent {\sc Proof of Proposition \ref{prop2.1}.} We first
rescale the immersion $\phi$ to $\tilde \phi=d_{2/r}\circ \phi$, where
$d_{2/r}$ is the dilation by the factor $2/r$. Thus the metric
induced by $\tilde x$ in $M$ is $d\tilde s^2=(4/r^2)ds^2$, where
$ds^2$ is the metric induced by $\phi$. We will denote the quantities
measured relative to the new metric by the superscript $\sim$.
Notice that the second fundamental form $\widetilde A$ satisfies
$|\widetilde A|^2=\frac{r^2}4\,|A|^2$.

Therefore, Proposition \ref{prop2.1} will be established once we
prove that given $\ve>0$ there exists $R_0$ such that, for
$r>R_0$,
$$
\sup_{M-\widetilde D_0(2)}|\widetilde A|^2<\ve/4.
$$

Given the above $\ve$, set $\ve_1<\ve/4$ and let $\delta>0$ be
given by Lemma \ref{lem2.4}. Since $M$ has strong finite total
curvature, there exists $R_0$ such that, for $r>R_0$,

$$
\begin{array}{rcl}
\delta &>&\left ( \int_{D_0(r/2,\infty)} |A|^q |\ro_0|^{q-n} \;dM\right)^{1/q}+\left ( \int_{D_0(r/2,\infty)} |\nabla |A||^q |\ro_0|^{2q-n} \;dM\right)^{1/q}\\[12pt]
&=&\left ( \int_{\widetilde D_0(1,\infty)} |\widetilde A|^q |\widetilde\ro_0|^{q-n} \;d\widetilde M\right)^{1/q}+\left ( \int_{\widetilde D_0(1,\infty)} |\widetilde\nabla |\widetilde A||^q |\widetilde\ro_0|^{2q-n} \;d\widetilde M\right)^{1/q}\hspace{-.5cm}\nonumber.
\end{array}
$$

For $x\in M-\widetilde D_0(2)$, we have $\widetilde
D_x(1)\subset\widetilde D_0(1,\infty)$ and then $\displaystyle{\inf_{\widetilde D_x (1)}\widetilde\ro_0 >1}$. Now, Lemma \ref{lem2.4}, with $\mu=|\widetilde\ro_0|^{q-n}$ and $\nu=|\widetilde\ro_0|^{2q-n}$, and the
above inequality imply that
$$
\sup_{\widetilde D_x(1/2)}|\widetilde A|^2<\ve_1,
$$
hence
$$
\sup_{M- \widetilde D_0(2)}|\widetilde A|^2\leq \ve_1 <\ve/4.
$$
This completes the proof of Proposition \ref{prop2.1}. \qed

\section{Uniqueness of the tangent plane at infinity}

The proof of  our Theorem \ref{thm1.1} depend on a series of
lemmas and a crucial proposition to be presented in a while.
 In this section, $\phi\colon M^n\to\re^{n+1}$ will always denote a complete hypersurface such that $\phi(M^n)$
passes through the origin $0$ of $\re^{n+1}$, with strong finite total curvature.

The following lemma is similar to Lemma \ref{lem2.3} in Anderson
\cite{A}.
\begin{lem}\label{lem3.1}
Let $\phi\colon M^n\to\re^{n+1}$ be as above and let
$r(p)=d(\phi(p),0)$, where $p\in M$ and $d$ is the distance in
$\re^{n+1}$. Then $\phi$ is proper and the gradient $\nabla r$ of $r$
in $M$ satisfies
$$
\lim_{r\to\infty}|\nabla r| = 1.
$$
In particular, there exists $r_0$ such that if $r>r_0$, $\nabla
r\ne0$, i.e., the function $r$ has no critical points outside the
ball $B(r_0)$.
\end{lem}

\begin{pf}
If the immersion is not proper, we can find a ray $\ga(s)$ issuing
from $0$ and parametrized by the arc length $s$ such that as $s$
goes to infinity the distance $r(\ga(s))$ is bounded. Let  such a
ray be given and set $T=\ga'(s)$. Let
$$
X=(1/2)\ov\nabla r^2=r\ov\nabla r,
$$
be the position vector field, where $\ov\nabla r$ is the gradient
of $r$ in $\re^{n+1}$. Then
$$
T\langle X,T\rangle=\langle\ov\nabla_TX,T\rangle+\langle
X,\ov\nabla_TT\rangle=1+\langle X,\ov\nabla_TT\rangle.
$$
Since $\ga$ is a geodesic in $M$, the tangent component of
$\ov\nabla_TT$ vanishes and
$$
\ov\nabla_TT=\langle\ov\nabla_TT,N\rangle
N=-\langle\ov\nabla_TN,T\rangle N=\langle A(T),T\rangle N.
$$
It follows, by Cauchy-Schwarz inequality, that
$$
|\langle X,\ov\nabla_TT\rangle|\le|X|\,|A(T)|\,|T|\le |X|\,|A|,
$$
hence
$$
T\langle X,T\rangle\ge1-|X|\,|A|.
$$

By using Proposition \ref{prop2.1} with $\ve=1/m^2$, and the facts
that $r=|X(s)|\le s$ and that $\ga$ is a minimizing geodesic, we obtain
\begin{equation}\label{eq3.1}
T\langle X,T\rangle(s)\ge 1-\frac1m,
\end{equation}
for all $s>R_0$, where $R_0$ is given by Proposition
\ref{prop2.1}. Integration of \eqref{eq3.1} from $R_0$ to $s$
gives
\begin{equation}\label{eq3.1.1}
\langle X,T\rangle(s)\ge\left(1-\frac1m\right)(s-R_0)+\langle
X,T\rangle(R_0).
\end{equation}

Because $r(s) = |X(s)| \ge \lgg X,T\rg(s)$, we see from (3.2) that
$r$ goes to infinity with $s$. This is a contradiction and proves
that $M$ is properly immersed.

Now let $\{p_i\}$ be a sequence of points in $M$ such that
$\{r(p_i)\} \to \infty$. Let $\ga_i$ be a minimizing geodesic from
$0$ to $p_i$, and denote again by $\ga(s)$ the ray which is the
limit of $\{\ga_i\}$. For each $\ga_i$, we apply the above
computation, and since
$$
\lgg X_i,T_i\rg(s) = \lgg r_i \ov\nabla r_i,T_i\rg(s) \le
r_i|\nabla r_i|(s),
$$
we have
$$
|\nabla r_i|(s) \ge \frac{\lgg X_i,T_i\rg(s)}{s} \ge \left(1-\frac
1m\right)\left(\frac{s-R_0}{s}\right) + \frac{\lgg
X_i,T_i\rg(R_0)}{s}\,,
$$
hence, for the ray $\ga(s)$,
\begin{equation}
|\nabla r|(s) \ge \left(1-\frac
1m\right)\left(\frac{s-R_0}{s}\right) + \frac{\lgg X,T\rg
(R_0)}{s}\cdot \tag{3.3}
\end{equation}

By taking the limit in (3.3) as $s \to \infty$, we obtain that
$\lim\limits_{s\to\infty} |\nabla r| \ge 1-\dfrac 1m\cdot$ Since
$m$ and the sequence $\{p_i\}$ are arbitrary, and $|\nabla r| \le
1$, we conclude that $\lim\limits_{r\to\infty} |\nabla r| = 1$,
and this completes the proof of Lemma 3.1. \qed
\end{pf}

\vglue .1in

\noindent{\bf Remark}. Related to Lemma 3.1,   Bessa, Jorge and
Montenegro \cite{BJM} proved independently that for
an immersion $\phi\colon M^n \to \re^N$ (of arbitrary codimension)
for which the norm $|\al|$ of the second fundamental form $\al$
satisfies
$$
\lim_{r\to\infty} \,\sup_{p\in M-D_0(r)} r^2|\al|^2 < 1
$$
it holds that $\phi$ is proper and that the distance function $r =
d(\phi(p),0)$, $p \in M$, has no critical point outside a certain
ball.

\medskip

Now, let $r_0$ be chosen so that $r$ has no critical points in
$W=\phi(M)-(B(r_0)\cap \phi(M))$. By Morse Theory, $x^{-1}(W)$  is
homeomorphic to $\phi^{-1}[\phi(M)\cap S(r_0)]\times[0,\infty)$. Let $V$
be a connected component of $\phi^{-1}(W)$, to be called an {\it
end\/} of $M$. It follows that $M$ has only a finite number of
ends. In what follows, we identity $V$ and $\phi(V)$.

Let $r>r_0$ and set
$$
\Sigma_r=\frac1r\,[V\cap S(r)]\subset S(1),
$$
$$
V_r=\frac1r\,[V\cap B(r)]\subset B(1).
$$
Denote by $A_r$ the second fundamental form of $V_r$. Then
$$
|A_r|^2(x)=r^2|A|^2(rx).
$$

\begin{lem}\label{lem3.2.2}
For $r>r_0$, $V\cap B(r)$ is connected.
\end{lem}

\begin{pf}
Notice that $V=S\times[0,\infty)$ where $S$ is a connected
component  of  $M\cap S(r_0)$. Assume that $V\cap B(r)$ has two
connected components, $V_1$ and $V_2$. Since $(V_1\cup V_2)\cap
S(r_0)$ is connected, either $V_1\cap S(r_0)$ or $V_2\cap S(r_0)$
is empty. Assume it is $V_2 \cap S(r_0)$.

Let $p\in V_2$. Since all the trajectories of $\nabla r$ start
from $V_1\cap S(r_0)$, there exists a trajectory $\vr(t)$ with
$\vr(0)\in V_1\cap S(r_0)$ and $\vr(t_2)=p$. Thus, there exist
$t_0,t_1\in [0,t_2]$, such that a trajectory of $\nabla r$
satisfies $|\vr(t_0)|=|\vr(t_1)|=r$. We claim that this implies
the existence of a critical point of $r$ at some point of
$\vr(t)$.

Indeed, let $f(t)=r(\vr(t))$. Then $f\colon\re\to\re$ is a smooth
function with $f(t_0)=f(t_1)$. Thus, there exists $\ov
t\in[t_0,t_1]$ with $f'(\ov t)=0$. But
$$
f'(t)=dr\left(\frac{d\vr}{dt} \right)= dr(\nabla r)=\langle\nabla
r, \nabla r\rangle.
$$
Therefore,
$$
0=f'(\ov t)=|\nabla r(\bar t)|^2
$$
and this proves our claim.

Thus we have reached a contradiction and this proves the lemma.
\end{pf}
\qed

\begin{lem}\label{lem3.3}
Let $0<\delta<1$ be given and fix a ring $A(\delta,1)\subset
B(1)$. Then, given $\ve>0$, there exists $r_1$ such that, for all
$r>r_1$ and all $x\in V_r\cap A(\delta,1)$, we have
$$
|A_r|^2(x)<\ve.
$$
\end{lem}

\begin{pf}
By Proposition \ref{prop2.1}, there exists $r_0$ such that for
$r>r_0$
\begin{equation}
r^2\sup_{x\in M-D_0(r)} |A|^2(x)<\delta^2 \ve. \tag{3.4}
\end{equation}
Take $r_1=r_0/\delta$. Then, for $r>r_1$ and $x\in V_r\cap A(\delta,1)$,
$$
r|x|>r\delta>r_0.
$$
Thus, by (3.4), for all $x\in V_r\cap A(\delta,1)$ and $r>r_1$,
\begin{equation}
r^2|x|^2\left[\sup_{y\in M-D_0(r|x|)}|A|^2(y)\right]<\delta^2\ve.
\tag{3.5}
\end{equation}
Now, by using again Proposition 2.2 and (3.5), we obtain that
$$
|A_r|^2(x)=r^2|A|^2(rx) \le r^2\sup_{y\in M-D_0(r|x|)}|A|^2(y) <
\frac{\delta^2\ve}{|x|^2}<\ve,
$$
for all $x \in V_r\cap A(\delta,1)$ and $r > r_1$, and this proves Lemma
\ref{lem3.3}. \qed
\end{pf}

By Lemma \ref{lem3.3}, we see that $|A_r|^2\to0$ uniformly in the
ring $A(\delta,1)$. It follows from this and the fact that $V_r$
is connected that we can apply Lemma 2.1(i) and conclude that a
subsequence $V_{r_i}$ of $V_r$, $r_i\to\infty$, converges $C^1$ to a union of
hypersurfaces $\pi$ in $A(\delta,1)$
 . Again, since $|A_r| \to 0$ uniformly, $\pi$ is a union of $n$-planes in $A(\delta,1)$ (see Remark after the proof of Lemma 2.1).
  Since $\delta$ is arbitrary, a subsequence again denoted by $V_{r_i}$ converges to $\pi$ in
  $B(1)-\{0\}$and the n-planes in $\pi$ all pass through the
  origin $0$. Thus, each two of them intersect along a linear
  (n-1)-subspace $L$ and the hypersurfaces $\Sigma_{r_i}\subset
  S^{n}(1)$, given by the inverse images of the regular values
  $r_i$ of the distance function $r$, converge to a family
  $\Sigma_\infty$ of equators of $S^{n}(1)$ each two of each
  intersect along $L\cup S^{n}(1)$. We claim that $\Sigma_\infty$
  contains only one equator. In fact, for $r_i$ large enough,
  by the basic
transversality theorem (\cite{Hi} Chapter 3, Theorem 2.1),
$\Sigma_{r_i}$ has a self intersection close
  to $L\cup S^{n}(1)$ and this contradicts the fact that
  $\Sigma_{r_i}$is an embedded hypersurface. It follows that $\pi$
  is a single n-plane passing through $0$, possibly with
  multiplicity $m\geq 1$. Since $\Sigma_\infty$ covers $S^{n-1}(1)$, which is
simply-connected, $m=1$. Thus $V$
  is embedded and $\pi$ is a single plane that
   passes through the origin.

The $n$-plane $\pi$ spanned by $\Sigma_\infty$ is called the {\it
tangent plane at infinity of the end $V$ associated to the
sequence\/} $\{r_i\}$. A crucial point in the proof of Theorem
\ref{thm1.1} is to show that this plane does not depend on the
sequence $\{r_i\}$. Here we use for the first time the hypothesis
on $H_n$.
\begin{prop}\label{lem3.4}
Each end $V$ of $M$ has a unique tangent plane at infinity.
\end{prop}

\begin{pf}
Suppose that $\{s_i\}$ and $\{r_i\}$, $s_i,r_i\to\infty$,  are
sequences of real numbers and that $\pi_1$ and $\pi_2$ are
distinct tangent planes at infinity associated to $\{s_i\}$ and
$\{r_i\}$, respectively. We can assume that the sequences satisfy
$$
s_1<r_1<s_2<r_2<\dots<s_i<r_i<\dots .
$$
Let $K$ be the closure of $B(3/4)-B(1/4)$ and let $N_1$ be the
normal to $\pi_1$, obtained as the limit of the normals to
$$
K\cap\left\{\frac1{s_i}V\right\}=\frac1{s_i}(V\cap s_iK).
$$
Similarly, let $N_2$ be the normal to $\pi_2$ obtained as the
limit of the normals to $K\cap\{(1/r_i)V\}$.

Now let $U_1$ and $U_2$ be neighborhoods in $S^n(1)$ of $N_1$ and
$N_2$, respectively, such that $U_1\cap U_2=\emptyset$. Thus,
there exists an index $i_0$ such that, for $i>i_0$, the normals to
$K_i^1=(s_iK)\cap V$ are in $U_1$ and the normals to
$K_i^2=(r_iK)\cap V$ are in $U_2$. If $K_i^1\cap
K_i^2\ne\emptyset$, for some $i>i_0$, this contradicts the fact
that $U_1\cap U_2=\emptyset$, and the proposition is proved.

Thus we may assume that, for all $i>i_0$, $K_i^1\cap
K_i^2=\emptyset$. In this case, we have $(1/4)r_i>(3/4)s_i$; here,
and in what follows, we always assume $i>i_0$. Set
$$
W_i=V\cap\left(B\left(\frac14 r_i\right)-B\left(\frac34
s_i\right)\right).
$$
Since $H_n$ does not change sign in $V$, we have that (\cite{HN},
Theorem II) $g(\po W_i)\supset\po(g(W_i))$. Since
$$
g\left(S\left(\frac 14 r_i\right)\cap V\right)\subset U_2,
$$
$$
g\left(S\left(\frac34 s_i\right)\cap V\right)\subset U_1,
$$
we have $g(\po W_i)\subset U_1\cup U_2$. Thus
\begin{equation}
\po(g(W_i))\subset g(\po W_i)\subset U_1\cup U_2\,. \tag{3.6}
\end{equation}

We claim that there exists a point $x\in \text{Int}(W_i)$ with
$H_n(x)\ne 0$.  Suppose that
\begin{equation}
\{x \in \text{ Int } W_i\,; H_n(x) \ne 0\} = \emptyset. \tag{3.7}
\end{equation}
Since $g(W_i)$ is connected and has nonvoid intersection with
$U_1$ and $U_2$ which are disjoint, there is a point $x_0 \in
\text{ Int }W_i$ such that $g(x_0) \notin U_1 \cup U_2$\,. Let
rank $A(x_0)=m$. By (3.7), $m < n$. Since the $k_i$'s are
continuous, there is a neighborhood $V$ of $x_0$ such that if $x
\in V$,
$$
n > \text{ rank } A(x) \ge m,
$$
where the left hand inequality follows from (3.7). This implies
that either rank $A$ is constant and equal to $m$ in a
neighborhood of $x_0$ or in each neighborhood of $x_0$ there is a
point such that the rank of $A$ at this point is greater than $m$.
In view of (3.7), the latter implies that we can find such a
point, to be called $y_0$, so that about $y_0$ there is a
neighborhood with rank $A = m_0 > m$.

In both cases, we obtain a point and a neighborhood of this point
for which rank $A$ is constant. Without loss of generality, we can
assume this point to be $y_0$. Notice that we can assume $g(y_0)
\notin U_1 \cup U_2$\,. By the Lemma of Chern-Lashof (\cite{CL},
Lemma 2), there passes through $y_0$ a piece $L^p$ of a
$p$-dimensional plane, $p = n-m_0$, along which $g$ is constant.
If $L^p$ intersects $\po W_i$, $g(y_0) \in g(\po W_i) \subset U_1
\cup U_2$, and this contradicts the choice of $y_0$. If not, a
point $\bar y_0$ in $\po L^p$ has again rank $A=m_0$ (\cite{CL},
Lemma 2), and arbitrarily close to $\bar y_0$, we have a point
$y_1$ and a neighborhood of $y_1$ whose rank is $m_1 > m_0$\,.
Thus, we can repeat the process.

After a finite number of steps, the process will lead either to
finding a point with rank $A=n$, what contradicts (3.7), or to
finding a piece $L$ of a plane of appropriate dimension with the
property that $L \cap \po W_i \ne \emptyset$. As we have seen
above, this is again a contradiction and proves our claim.

Thus, we can assume that there is a point $x \in \text{ Int}(W_i)$
with $H_n(x) \ne 0$. Then $g(W_i)$ contains an open
 set around $g(x)$. We can assume that $U_1$ and $U_2$ are small enough so that $g(x) \notin U_1 \cup U_2$. Since
 $g(W_i)$ is connected and has nonvoid intersection with both $U_1$ and $U_2$, the fact that there are interior points in  $g(W_i)$ and (3.6)
  imply that
\begin{equation}
g(W_i) \supset S^n(1) - \{U_1 \cup U_2\}. \tag{3.8}
\end{equation}

On the other hand, because
$$
\big(\Sigma k_i^2\big)^q > Ck_1^2 \dots k_n^2,
$$
for a constant $C=C(n)$, we have that
$$
|H_n|< \frac1{\sqrt C}\,|A|^q.
$$
Furthermore, since $\phi$ has strong finite total curvature,
$$
\int_{W_i}|A|^q |\rho_0|^{n-q}\, dM\to0,\quad i\to\infty.
$$

Therefore, since
$$
\text{Area } g(W_i)\le\int_{W_i}|H_n|\, dM < (\frac1{\sqrt
C})\int_{W_i}|A|^q |\rho_0|^{n-q}\, dM,
$$
we have that Area $g(W_i)\to0$. This a contradiction to (3.8), and
completes the proof of  Proposition 3.4. \qed
\end{pf}

\section{Proofs of Theorems \ref{thm1.1}, \ref{thm1.2} and \ref{thm1.3}}

\noindent {\sc Proof of Theorem \ref{thm1.1}.}
  (i) has already been proved in Lemma 3.1. To prove (ii), we apply to each end $V_i$ the inversion $I\colon \re^{n+1}-\{0\}\to\re^{n+1}-\{0\}$, $I(x)=x/|x|^2$. Then  $I(V_i)\subset B(1)-B(0)$ and as $|x|\to\infty$ in $V_i$, $I(x)$ converges to the origin $0$. It follows that each $V_i$ can be compactified with a point $q_i$. Doing this for each $V_i$, we obtain a compact manifold $\ov M$ such that $\ov M-\{q_1,\dots,q_k\}$ is diffeomorphic to $M$. This prove (ii).

To prove (iii), we use again the above inversion and observe that,
by Proposition 3.4, as $|x|\to \infty$ in $V_i$, the normals at
$I(x)$ converge to a unique normal $p_i\in S_1^n$ (namely, to the
normal of the  unique plane at infinity of $V_i$). Thus we obtain
a continuous extension $\ov g\colon\ov M\to S_1^n$ of $g$ by
setting $\ov g(q_i)=p_i$. This proves (iii).

\vskip .1in

\noindent {\sc Proof of Theorem 1.2{\rm(i)}}. We first observe
that $S_1^n-(N)$ is still simply-connected. This comes from the
fact that a closed curve $C$ in $S_1^n-(N)$ is homotopic to a
simple one and a disk generated by such a curve can, by
transversality, be made disjoint of $N$ by a small perturbation.
Thus $C$ is homotopic to a point in $S_1^n-(N)$.

Next, the restriction map
$$
\tilde g\colon M-\overline g^{-1}(N\cup\{p_i\}) \to S_1^n - (N\cup
\{p_i\})
$$
where $p_i$ is defined in the proof of Theorem 1.1, is clearly
proper and its Jacobian never vanishes. In this situation, it is
known that the map is surjective and a covering map (\cite{WG},
Corollary 1). Since $S_1^n-(N\cup \{p_i\})$ is simply-connected,
$\tilde g$ is a global diffeomorphism.

To complete the proof we must show that if $\ov g(n_1)=\ov
g(n_2)=p$, $n_1,n_2\in \overline g^{-1} (N\cup\{p_i\})$ then
$n_1=n_2$. Suppose that $n_1\ne n_2$. Let $W\subset S^n(1)$ be a
neighborhood of $p$. By continuity, there exist disjoint
neighborhoods $U_1$ of $n_1$ and $U_2$ of $n_2$ in $\ov M$ such
that $\ov g(U_1)\subset W$ and $\ov g(U_2)\subset W$. Choose $t\in
\ov g(U_1)\cap\ov g(U_2)$, $t\notin N\cup\{p_i\}$. Then, there
exist $r_1\in U_1$ and $r_2\in U_2$ such that $\tilde
g(r_1)=\tilde g(r_2)=t$. But this contradicts the fact that
$\tilde g$ is a diffeomorphism and concludes the proof of (i).

\medskip

(ii) We will use a result of Barbosa, Fukuoka and Mercuri
\cite{BFM}. By using Hopf's theorem that the Euler characteristic
$\chi(\ov M)$ of $\ov M$ is equal to the sum of the indices of a
vector field, the following expression is obtained in \cite{BFM}
Theorem 2.3: if $n$ is even,
$$
\chi(\ov M) = \sum_{i=1}^k (1+I(q_i)) + 2d\sigma.
$$
Here $I(q_i)$ is the multiplicity of the end $V_i$ (since $n \ge
3$, $I(q_i)=1$ in our case), $\sigma$ is $\pm1$ depending on the
sign of $H_n$\,, $k$ is the number of ends and $d$ is the degree
of the Gauss map $\overline g$. From Theorem 1.2 (i), $\overline
g$ is a homeomorphism. Thus, $d=1$ and, since $n$ is even,
$\chi(\ov M) = 2$. It follows that
$$
2 = 2k + 2\sigma.
$$
Thus $k=2$ and $\sigma = -1$, and the result follows. \qed

\medskip

\noindent{\sc Proof of the Gap Theorem}. First, we easily compute
that
$$
|A|^{2n} > (n!)H_n^2.
$$
Thus, since $H_n$ is the determinant of the Gauss map $g\colon M^n
\to S_1^n$\,, we obtain
$$
\int_M|A|^n\,dM > \sqrt{n!} \int_M |H_n|dM = \sqrt{n!} \text{ area
of $g(M)$ with multiplicity.}
$$

The extended map $\ov g\colon \ov M \to S_1^n$\,, which is given
by Theorem 1.1, has a well defined degree $d$, hence
$$
\text{area } g(M) = \text{ area } \ov g(\ov M) = d \text{ area }
S_1^n\,.
$$
Now, assume that $\phi(M)$ is not a hyperplane. We claim that
$d\ne0$. To see that, we first show that there exists a point in
$M$ where $H_n \ne 0$.

Suppose the contrary holds. Then, since $\phi(M)$ is not a
hyperplane, there is a point $x_\ell \in M$ such that rank $A$ at
$x_\ell$ is $\ell$, $0 < \ell < n$. Thus, by using the Lemma of
Chern-Lashof (\cite{CL}, Lemma 2) in the same way as we did in
Proposition 3.4, we arrive, after a finite number of steps, at one
of the two following situations. Either we find a point where $H_n
\ne 0$, which is a contradiction, or we find an open set $U_j
\subset M$, whose points satisfy rank $A = j \ge \ell$, $j < n$,
foliated by $(n-j)$-planes the leaves of which extend to infinity.
In the second situation, observe that the Gauss map on each leaf
is constant and, since there is only one normal at infinity for
each end, the normal map is constant on $U_j$\,. Thus $U_j$ is a
piece of a hyperplane, and we find again a contradiction, this
time to the fact that $n > j \ge \ell > 0$.

Therefore, there exists a point $x_0 \in M$ with $H_n(x_0) \ne 0$.
Then, for a neighborhood $V$ of $x_0$\,, we have that $H_n(x) \ne
0$, $x \in V$, and that $g(V) \subset S_1^n$ is a neighborhood of
$g(x_0)$. By Sard's theorem, the set of critical values of $g$ has
measure zero, hence some point of $g(V)$ is a regular value. It
follows that the Gauss map $g$ has regular values whose inverse
images are not empty. Since $H_n$ does not change sign, this prove
our claim.

Furthermore the area $\sigma_n$ of a unit sphere of $\re^{n+1}$ is
given by
$$
\sigma_n = \frac{2(\sqrt \pi)^{n+1}}{\Ga((n+1)/2)}\,;
$$
here $\Ga$ is the gamma function, which, in the present case is
given by
\begin{align*}
&\Ga((n+1)/2) = ((n-1)/2)!, \text{ if $n$ is odd}\\
&\Ga((n+1)/2) = \frac{(n-1)(n-3)\dots 1}{2^{n/2}}\,\sqrt\pi,
\text{ if $n$ is even.}
\end{align*}

It follows that, for all non-planar $x \in C^n$,
$$
\int_M |A|^n\,dM > 2\sqrt{n!}\, (\sqrt\pi)^{n+1}\big/\Ga((n+1)/2).
\qquad\quad\qed
$$

{\bf Acknowledgement:} We thank Elon Lima, Barbara Nelli, Wladimir Neves, Harold Rosenberg and Walcy Santos, Frank Pacard and Heudson Mirandola for conversations related to this paper.

\vskip .5in

\begin{center}
\begin{tabular}{l l}
M. do Carmo \qquad & \qquad M.F. Elbert \\
 IMPA  \qquad & \qquad  UFRJ \\
Estrada Dona Castorina 110 \qquad & \qquad  Instituto de Matem\'atica \\
 CEP 22460-320  \qquad & \qquad Cx. Postal 68530, CEP 21941-909  \\
 Rio de Janeiro, RJ, Brasil \qquad & \qquad Rio de Janeiro, RJ, Brasil \\
e-mail: manfredo@impa.br \qquad & \qquad e-mail: fernanda@im.ufrj.br\\
\end{tabular}
\end{center}

\end{document}